\documentclass{amsart}
\usepackage{amssymb,latexsym}

\theoremstyle{plain}
\newtheorem{theorem}{Theorem}

\newtheorem{lemma}{Lemma}
\newtheorem{proposition}{proposition}
\newtheorem{remark}{Remark}

\newtheorem{example}[theorem]{Example}

\theoremstyle{definition}
\newtheorem{definition}{Definition}

\theoremstyle{remark}

\numberwithin{equation}{section}

\begin{document}
\title[Fuzzy signed measure ]{Fuzzy Signed Measure}
\author{Jun Tanaka}
\address{University of California, Riverside, USA}
\email{juntanaka@math.ucr.edu, yonigeninnin@gmail.com, junextension@hotmail.com}

\date{January, 10, 2008}

\maketitle

\section{Introduction}
It is well known that one could obtain the Caratheodory Extension Theorem on fuzzy measurable space. Now a fuzzy measure is defined on $\sigma$-algebras. As one already saw, Fuzzy measure is a classical measure, provided that fuzzy sets are restricted to classical sets. As the classical measure theory goes, we will define a fuzzy signed measure on $\sigma$-algebras, as well as positive and negative sets. Herein, we will show that the Fuzzy Hahn Decomposition Theorem, which is a generalization of the classical Hahn Decomposition Theorem, decompose any space X into a positive set A and a negative set B such that A+B=X and the signed measure of $A \wedge B $ is 0.

\section{Preliminaries}

In this section, we shall briefly review the well know facts about lattice theory (e.g. Birkhoff [1], Iwamura), propose an extension lattice, and investigate its properties. Later in this section, we shall review Caratheodory Extension Theorem on fuzzy measurable sets. (L,$\wedge  $,$\vee  $) or simply L, under closed operations $\wedge  $ $\vee  $, is called a lattice. If it satisfies in addition to the distributive law is called a lattice. For two lattices L and L', a bijection from L to L', which preserves lattice operations are called a lattice isomorphism, or simply an isomorphism. If there is an isomorphism from L to L', then L is called lattice-isomorphic with L', and we write L $\cong$ L'. We write x $\leq $ y if x $\wedge$ y = x or, equivalently, if x $\wedge$ y = y. L is called complete, if any subset A of L includes the supremum $\vee $A, and infimum $\wedge$A, with respect to the above order. A complete lattice L includes the maximum and minimum elements, which are denoted by I and O, or 1 and 0, respectively [1].

\begin{definition}

Unless otherwise stated, X is a space and $\mu_{\Box} $ is a membership function of any fuzzy set $\Box$. If a family $\sigma$ of membership functions on X satisfies the following conditions, then it is called fuzzy $\sigma$-algebra;

(1) $\forall  \alpha \in $ [0,1], $ \alpha $ is constant;   $ \alpha \in \sigma$.

(2)  $\forall  \mu \in \sigma$, 1-$ \mu \in \sigma$.

(3) if $(\mu_{n} )_{n \in N}  \in \sigma^{N}$, then $ \sup \mu_{n} \in \sigma  $.

\end{definition}

\begin{definition}\label{De:2}
If m : $\sigma$ $\mapsto$ $\mathbb{R} \cup \{\infty \}$ satisfying the following properties, then m is called a fuzzy measure.

(1) m $ (\emptyset) $ = m $ (0)$ = 0.

(2) $\forall \mu , \eta  \in \sigma$ s.t. $m( \mu) , m(\eta) \geq   $ 0 : $\mu \leq \eta$ $ \Rightarrow  $  m $ ( \mu) \leq $ m $(\eta)$.

(3) $\forall \mu , \eta  \in \sigma$ : m$( \mu \vee \eta ) + $ m$( \mu \wedge \eta ) = $ m$( \mu)  +  $ m$ (\eta)   $.

(4) $( \mu_{n} )_{n \in N} $ $\subset  \sigma^{N} $ such that $\mu_{1} \leq   \mu_{2} \leq    \cdots  \leq  \mu_{n} \leq  \cdots $ : $\sup \mu_{n} = \mu$ $\Rightarrow$ m$ ( \mu ) =  \lim $ m$( \mu_{n} )$.

\end{definition}

\begin{definition}\label{De:3}
By an outer fuzzy measure $m^{\ast}$, we mean an extended real-value set function defined on [0,1]$^{X}$, having the following properties:

(1) $m^{\ast}$(0) = 0

(2) $m^{\ast}$  $ ( \mu) \leq $  $m^{\ast}$ $(\eta)$ for $\mu \leq \eta$

(3) $m^{\ast} ( \vee_{i=1}^{\infty}  \mu_{E_{i}} )   $   $ \leq $ $\vee_{i=1}^{\infty}   m^{\ast}(  \mu_{E_{i}} )$.

\end{definition}

\begin{example}
Suppose
\[
m^{\ast}(\mu_{E} ) =
\begin{cases}
0  \quad &  \mu_{E} = 0 \\
1  \quad &  \mu_{E}  \neq  0
\end{cases}
\]

Then $m^{\ast}$ is an outer fuzzy measure which is not fuzzy measure on [0,1]$^X$, if X has at least two points.

\end{example}

\begin{proposition}\label{Pro:1}
 Let F be a class of fuzzy subsets of X containing empty set such that for every $\mu_{A} \leq \mu_{X}$, there exists a sequence $(\mu_{B_{n}})^{\infty}_{n=1} $, $\mu_{B_{n}} \in F$ such that $\mu_{A} \leq$ $\sup(\mu_{B_{n}})^{\infty}_{n=1}$. Let $\tau$ be an extended real-valued function on F such that $\tau (0)$ = 0 and $\tau(\mu_{A}) \geq$ 0 for $\mu_{A} \in $ F. Then $m^{\ast}$ defined on [0,1]$^{X}$ by $m^{\ast} (\mu_{A})$ = $\inf \{ \tau ( \vee^{\infty}_{i=1} \mu_{B_{n}}) : \mu_{B_{n}} \in F $ s.t. $ \mu_{A} \leq   \vee^{\infty}_{i=1} \mu_{B_{n}}      \}  $ is a outer fuzzy measure.

\begin{proof}
Clearly, $m^{\ast} ( \emptyset)$ = 0. Secondly, if $\mu_{A_{1}} \leq \mu_{A_{2}}$ and $\mu_{A_{2}}  \leq  \vee^{\infty}_{i=1} \mu_{B_{n}} $, then $\mu_{A_{1}}  \leq  \vee^{\infty}_{i=1} \mu_{B_{n}} $. Thus $m^{\ast} (\mu_{A_{1}}) \leq m^{\ast} (\mu_{A_{2}})$.

Finally, let $\mu_{E_{n}} \leq \mu_{X}$ for each natural number n. Then $m^{\ast}$($\mu_{E_{n}}$) = $\infty $ for some n, $m^{\ast}$($\vee_{n=1}^{\infty}  \mu_{E_{n}}$) $\leq$ $\vee_{n=1}^{\infty} m^{\ast} \mu_{E_{n}}  $.

\end{proof}
\end{proposition}

The following theorem is an extension of the above proposition.

\begin{theorem}\label{Theo:1}
The class B of $m^{\ast}$-fuzzy measurable sets is a $\sigma$-algebra. Also $\overline{m}$, the restriction of $m^{\ast}$ to B, is a measure.

\end{theorem}

\begin{theorem}\label{Theo:2}
Let m be a fuzzy measure on a $\sigma$-algebra $\sigma \subseteq $ [0,1]$^X$.

Suppose for $\mu_{E} \leq \mu_{X}$
\[ m^{\ast} ( \mu_{E} )  = \inf \{ m ( \vee^{\infty}_{i=1} \mu_{E_{n}}) : \mu_{E_{n}} \in \sigma  s.t.  \mu_{E} \leq   \vee^{\infty}_{i=1} \mu_{E_{n}}      \}
\]

Then the following properties hold:

(i) $ m^{\ast}$ is an outer fuzzy measure.

(ii) $ \mu_{E} \in \sigma$ implies m($\mu_{E}  $) = $ m^{\ast}$($\mu_{E}  $)

(iii) $ \mu_{E} \in \sigma$ implies  $ \mu_{E}$ is $m^{\ast}$-fuzzy measurable.

(iv) The restriction $\overline{m}$ of $m^{\ast}$ to the $m^{\ast}$-fuzzy measurable sets in an extension of m to a fuzzy measure on a fuzzy $\sigma$-algebra containing $\sigma$.

(v) If m is fuzzy $\sigma$-finite, then $\overline{m}$ is the only fuzzy measure (on the smallest fuzzy $\sigma$-algebra containing $\sigma$) that is an extension of m.

\end{theorem}

\section{Fuzzy Signed Measure}

Let $m_{1}$ and $m_{2}$ be fuzzy measures defined on the same $\sigma$-algebra $\sigma$. If one of them is finite, the set function $m ( \mu_{E} ) =  m_{1}( \mu_{E} )  - m_{2}( \mu_{E} )$ , $\mu_{E} \in \sigma $ is well defined and countably additive on $\sigma$. However, it is not necessarily nonnegative; it is called a signed measure.

\begin{definition}\label{De:1}
 By a fuzzy signed measure on the fuzzy measurable space (X, $\sigma$) we mean $\nu$ : $\sigma$ $\mapsto$ $\mathbb{R} \cup \{\infty \}$ or $\mathbb{R} \cup \{- \infty \}$, satisfying the following property:

(1) $\nu (\emptyset) $ = $\nu (0)$ = 0.

(2) $\forall \mu , \eta  \in \sigma$ s.t. $\nu( \mu) , \nu(\eta) \geq   $ 0 : $\mu \leq \eta$ $ \Rightarrow  $  $ \nu( \mu) \leq \nu(\eta)$.

\ \ \ \ \ $\forall \mu , \eta  \in \sigma$ s.t. $\nu( \mu) , \nu(\eta) \leq   $ 0 : $\mu \leq \eta$ $ \Rightarrow  $  $  \nu(\eta)  \leq  \nu( \mu)$.

(3) $\forall \mu , \eta  \in \sigma$ : $\nu( \mu \vee \eta ) +   \nu( \mu \wedge \eta ) = \nu( \mu)  +   \nu(\eta)   $.

(4) $( \mu_{n} )_{n \in N} $ $\subset  \sigma^{N} $ such that $\mu_{1} \leq   \mu_{2} \leq    \cdots  \leq  \mu_{n} \leq  \cdots $ : $\sup \mu_{n} = \mu$ $\Rightarrow$ $ \nu( \mu ) =  \lim \nu( \mu_{n} )$.

This is meant in the sense that if the left-hand side is finite, the series on the right-hand side is convergent, and if the left-hand side is $ \pm \infty $, then the series on the right-hand side diverges accordingly.
\end{definition}

\begin{remark}
The Fuzzy signed measure is a fuzzy measure when it takes only positive value. Thus, the fuzzy signed measure is a generalization of fuzzy measure.

\end{remark}

\begin{definition}
A is a positive fuzzy set if for any fuzzy measurable set E in A, $\nu (\mu_{E}) \geq 0$. Similarly, B is a negative fuzzy set if for any fuzzy measurable set E in B, $\nu (\mu_{E}) \leq 0$.

\end{definition}

\begin{lemma}\label{Le:1}
Every fuzzy subset of a positive fuzzy set is a positive fuzzy set and any countable union of positive fuzzy sets is a positive fuzzy set.

\begin{proof}
The first claim is clear. Before we show the second claim, we need to show that every union of positive sets is a positive set. Let A, B be fuzzy positive sets and E $\leq$ A $\vee $ B be a fuzzy measurable set. By (2) in Definition $\ref{De:1}$, 0 $\leq  $ $\nu (\mu_{B} \wedge \mu_{E}) $ - $\nu (\mu_{A} \wedge \mu_{B} \wedge \mu_{E}) $. By (3), $\nu (\mu_{E}) \geq$ 0. Now by induction, every finite union of fuzzy positive sets is a positive fuzzy set. Let $A_{n}$ be a positive fuzzy set for all n and E  $\leq$ $\vee $ $A_{n}$ be a fuzzy measurable set. Then $\mu_{E_{m}}$ := $ \mu_{E} \wedge \vee_{n=1}^{m} \mu_{A_{n} } $ =  $\vee_{n=1}^{m} \mu_{ E}  \wedge \mu_{A_{n} } $. Then $E_{m}$ is a fuzzy measurable set and a positive fuzzy set. In particular,  $\mu_{E_{m}} \leq \mu_{E_{m+1}}$ for all n and $ \mu_{E}$ = $ \vee_{n=1}^{\infty} \mu_{E_{m} }$. Thus 0 $\leq  \lim \nu ( \mu_{E_{m}} ) $ = $\nu ( \mu_{E }) $. Therefore $\vee $ $A_{n}$ is a positive set.

\end{proof}
\end{lemma}

\begin{lemma}\label{Le:2}

Let E be a fuzzy measurable set such that 0 $< \nu (\mu_{E})   <  \infty$. Then there is a positive fuzzy set A $\leq E$ with $ \nu (\mu_{A}) >$ 0.

\begin{proof}
If E is a positive fuzzy set, we take A=E. Otherwise, E contains a set of negative measure. Let $n_{1}$ be the smallest positive integer such that there is a measurable set $E_{1}  \subset $ E with $ \nu ( \mu_{E_{1} }) < - \frac{1}{n_{1}} $. Proceeding inductively, if E$\wedge \wedge_{j=1}^{k-1} E_{j}^C  $ is not already a positive set, let $n_{k}$ be the smallest positive integer for which there is a fuzzy measurable set $ E_{k} $ such that $ E_{k} \leq E \wedge \wedge_{j=1}^{k-1} E_{j} $ and $\nu ( \mu_{E_{k}}) < - \frac{1}{n_{k}} $.

Let A = $( \vee  E_{k})^C $.

Then $\nu (\mu_{E})$ = $ \nu ( \mu_{E} \wedge \mu_{A}  ) + \nu(\mu_{E}   \wedge   \vee  \mu_{E_{k} } )  $ = $ \nu (\mu_{E} \wedge \mu_{ A } ) + \nu(   \vee  \mu_{E_{k}}  )  $. Since  $\nu (\mu_{E})$ is finite, $ \lim_{n \rightarrow \infty}   \nu(   \vee^{n} \mu_{E_{k}}  )$ is finite and $ \nu(   \vee \mu_{ E_{k} } ) \leq$ 0. Since $\nu (\mu_{E}) >$ 0 and $ \nu(   \vee  \mu_{E_{k}}  ) \leq$ 0, $ \nu (\mu_{E} \wedge \mu_{A}  ) >$ 0.

We will show that A is a positive set. Let $  \epsilon  >$ 0. Since $ \frac{1}{n_{k}} \rightarrow  $ 0, we may choose k such that $-\frac{1}{n_{k} - 1}    $, which is greater than $ - \epsilon $. Thus A contains no fuzzy measurable sets of measure less than $ - \epsilon $. Since $ \epsilon $ was arbitrary positive number, it follows that A can contain no sets of negative measure and so must be a positive fuzzy set.

\end{proof}
\end{lemma}

\section{Fuzzy Hahn Decomposition}

Without loss of generality, let's omit + $\infty $ value of $\nu$. Let $\lambda$ = $\sup \{ \nu (\mu_{A}) : A $ is a fuzzy positive set $      \}   $.

Then $\lambda   \geq$ 0 since $\nu (\emptyset) $ = 0.

Let {$A_{i} $} be a sequence of positive fuzzy sets such that $\lambda$ = $\lim  \nu (\mu_{A_{i}} )$ and $\mu_{  A}$ = $\vee \mu_{A_{i} }$. By Theorem $\ref{Le:1}$, A is a positive fuzzy set and $\lambda \geq \nu (\mu_{A})$.

$\vee^{n} \mu_{ A_{i}} \leq $ $\mu_{A}$ for any n implies  $\nu ( \vee^{n} \mu_{ A_{i}}) \geq$ 0 for any n. Thus $\lambda$ = $\lim \nu (  \mu_{ A_{i} })$ = $ \nu (\mu_{A}) $ = 0.

Let $E \leq A^C$ be a positive set. Then $\nu (\mu_{E}) \geq$ 0 and $A \vee E $ is a positive fuzzy set. Thus $\lambda   \geq$ $ \nu ( \mu_{A} \vee \mu_{E}    ) = \nu ( \mu_{A} ) + \nu ( \mu_{E} )    - \nu (\mu_{ A } \wedge \mu_{E}  )   $ = $\lambda + \nu ( \mu_{E}) - \nu ( \mu_{A } \wedge \mu_{E}  )    $. Thus $\nu (\mu_{E}) = \nu ( \mu_{A}  \wedge \mu_{E}  )    $. $\nu (\mu_{E})$ = 0 since $ \mu_{A } \wedge \mu_{E} \leq \mu_{A} \wedge \mu_{A}^C$ and $\nu ( \mu_{A} \wedge \mu_{A}^C ) $ = 0.

Thus, $A^C$ contains no fuzzy positive subsets of positive measure and hence no subsets of positive measure by Lemma $\ref{Le:2}$. Consequently, $A^C$ is a negative fuzzy set.

\section{conclusion}
Let X be a space. Then by the previous theorem, we find such a positive fuzzy set A and a negative fuzzy set B (= $A^C$). By the fuzzy measurability of $\nu$, $\nu ( \mu_{A}  \wedge   \mu_{A}^C    ) $ = 0. $A + A^C$ = X in the sense that $\mu_{A} + \mu_{A^C}$ = 1. These characteristics provides X = $ A \cup  B$ and $ A \cap B = \emptyset$ in the classical set sense.



\end{document}